\documentclass[12pt,a4paper]{article}
\usepackage[utf8]{inputenc}
\usepackage[T1]{fontenc}
\usepackage[english]{babel}
\usepackage{amsmath}
\usepackage{amsfonts}
\usepackage{amssymb}
\usepackage{graphicx}
\usepackage{amsthm}
\usepackage{hyperref}
\usepackage[open,openlevel=1]{bookmark}
\newtheorem{theorem}{Theorem}[section]

\newtheorem{corollary}{Corollary}[section]

\newtheorem*{acknowledgement}{Acknowledgement}
\title{Two Applications of Topological Fixed Point Theory}

\author{Chaitanya Ambi\footnote{Chennai Mathematical Institute, H1, SIPCOT IT Park, Siruseri, Kelambakkam 603103, India. \newline Email: chaitanya.ambi@gmail.com}}
\begin{document}
 \maketitle
\begin{abstract}
\indent

We give a new proof of Cartan's fixed point theorem using topological fixed point theory. For an odd dimensional, simply connected and complete manifold having non-positive curvature, we further prove that every isometry with finite order not only fixes a point, but also maps a \textit{pencil} of geodesics to itself.\\

We apply similar techniques to the automorphism group of a smoothly bounded domain in the complex plane having connectivity three or more. We show that any prime dividing the order of this finite group must divide a certain integer depending only on the connectivity of the domain.
\end{abstract}
\section{Introduction.}
\indent

    The geometric properties of a complete Riemannian manifold $M$ with non-positive curvature are intimately linked with its fundamental group. Specifically, the structure of the isometry group $Isom(M)$ has been a field of intense investigation. It is not easy to ascertain whether a geometric property is, in fact, topological. As an attempt in this direction, we present a new proof of the following version of Cartan's fixed point theorem (see \cite{petersen}, Chap. 6, Thm. 23, p. 164). A stronger version of the same for precompact subgroups of $Isom(M)$  is known (see \cite{eberlein}, Thm. 1.4.6, p. 21).
\begin{theorem}[E. Cartan]\label{thm-Cartan}
 Consider a complete, simply connected Riemannian manifold $M$ of non-positive curvature. If $\sigma \in Isom(M)$ is an element of finite order, then \begin{equation*}
                                                                                                                                                       Fix(\sigma) \not = \emptyset .
                                                                                                                                                      \end{equation*}

\end{theorem}
\indent

    The above theorem can be proved by considering the centre of mass of the orbit of $\sigma$. In this article, we present a proof using topological fixed point theory. Our approach also reveals the following fact which is interesting by its own. By a \textit{pencil} of geodesics, we shall mean a (non-intersecting) family of geodesics any two of which lie within a finite Hausdorff distance of one another.
\begin{corollary}\label{cor-Cartan}
 In addition to the hypotheses of Cartan's Theorem \ref{thm-Cartan}, assume that $\dim M$ is odd. Then there exists a pencil $\Gamma$ of geodesics such that \begin{equation*}
                                                                                                                                                              \sigma(\Gamma)= \Gamma.
                                                                                                                                                             \end{equation*}

\end{corollary}
\indent

    Next, certain domains in the complex space $\mathbb{C}^{n}$ also become complete manifolds with non-positive curvature when equipped with the Bergmann metric. In this vein, we consider the group of conformal automorphisms of a bounded planar domain (with smooth boundary) having finite connectivity (see \cite{krantz}, Def. 4.2.2, p. 85). It is well-known that such a group is finite whenever the connectivity exceeds two (see \cite{krantz}, Thm. 12.2.3, p. 263). We prove the following result about such domains (where $Aut(.)$ denotes the group of conformal automorphisms of the respective domain).
\begin{theorem}\label{thm-auto}
 Let $\Omega_{k} \subset \mathbb{C}$ be any bounded planar domain having connectivity $k+1$ (with $k\geq 2$) whose boundary consists of finitely many smooth Jordan curves. If $p$ is a prime integer dividing the order of $Aut(\Omega_{k})$, then we have \begin{equation*}
         k^{3} \equiv k \mod p .                                                                                                                                                       \end{equation*}
\end{theorem}
\indent

     As far as the author knows, it is an open problem to characterise the finite groups which may occur as the automorphism groups of planar domains. In particular, by \cite{heins}, we have the following sharp upper bound for $k\geq 2$.
 \begin{align*}
 \#Aut(\Omega_{k}) \leq 2k  \quad \text{ for } k\not = 4,6,8,12,20,\\
 \#Aut(\Omega_{4}) \leq 12,\\
 \#Aut(\Omega_{6}),\#Aut(\Omega_{8}) \leq  24,\\
 \#Aut(\Omega_{12}),\#Aut(\Omega_{20}) \leq 60.
\end{align*}
\indent

    Theorem \ref{thm-auto} shows that there are further restrictions on the possible order of an element in $Aut(\Omega_{k})$. Hence, the subgroups of $Aut(\Omega_{k})$ are constrained by the value of $k$. For instance, the cyclic group of order $7$ cannot occur as a subgroup of $Aut(\Omega_{k})$ for any $k < 6$ (otherwise, we would have $7 | (k^{3}-k)$). Note that the upper bound mentioned above does not exclude this possibility as $\#Aut(\Omega_{5})\leq 10$.
\section{Proof of Theorem \ref{thm-Cartan}.}
\indent

    By the Cartan-Hadamard Theorem (see \cite{petersen}, Chap. 6, Thm. 22, p. 162), $M$ is diffeomorphic to $\mathbb{R}^{n}$, where $n=\dim M$. We compactify $M$ by adding its ideal boundary $M(\infty)$ to obtain the \textit{Gromov compactification}
    \begin{equation*}
     \bar{M}:= M \cup M(\infty).
    \end{equation*}
    We refer the reader to \cite{gromov}, Chap, 3, pp. 22-32 for the details. In fact, $\bar{M}$ is homeomorphic to the closed unit ball $\mathbb{B}^{n}\subset \mathbb{R}^{n}$.\\

    Next, let $\sigma \in Isom(M)$ have a finite order. It is known that $\sigma$ induces a homeomorphism $\bar{\sigma}$ of $\bar{M}\approx \mathbb{B}^{n}$.\\

    We shall argue inductively to reduce the dimension $n$. By Brouwer's fixed point theorem, $\bar{\sigma}$ has a fixed point $x \in\bar{M}$. If $x \in M(\infty)$, we observe that (see \cite{gromov}, Chap. 3, p. 30) $\bar{\sigma}$ maps a \textit{horosphere} centred at $x$ to itself. This horosphere in $\bar{M}$ is diffeomorphic to $\mathbb{B}^{n-1}$. So, we may repeat the whole argument to keep reducing the dimension. At each stage, the restriction of $\bar{\sigma}$ fixes a point lying either in $M$ or on $M(\infty)$. \\

    If the latter case were to occur until the dimension reaches two, the argument would yield a geodesic $\gamma \subset M$ satisfying
    \begin{equation*}
     \sigma(\gamma)=\gamma.
    \end{equation*}
    But now, $\gamma$ is \textit{isometric} with $\mathbb{R}$ (equipped with the standard Euclidean metric). Any isometry of the latter either fixes a point or has infinite order. As $\sigma$ is of finite order, we get a fixed point as asserted. This completes the argument.
    \subsection{Proof of Corollary \ref{cor-Cartan}.}
\indent

    As in the proof of Theorem \ref{thm-Cartan}, the extended map $\bar{\sigma}$ restricts to a homeomorphism of $M(\infty)$. Since we want to find a family of undirected geodesics (and not geodesic rays), we may well disregard the directions to obtain the quotient of $M(\infty)$ by $\mathbb{Z}/2\mathbb{Z}$ action. As $M(\infty)\approx \mathbb{S}^{n-1} $, the quotient space is homeomorphic to the real projective space $\mathbb{R}P^{n-1}$. This space also has an action induced by that of $\bar{\sigma}$.\\

    Assume the contrary that the above action induced by $\bar{\sigma}$ has no fixed points. Let $m\in \mathbb{N}$ be the smallest power of $\bar{\sigma}$ whose Lefschetz number is nonzero. It is clear that $m$ cannot exceed the order of $\bar{\sigma}$ because $(n-1)$ is even and so
    \begin{equation*}
     \lambda(Id_{\mathbb{R}P^{n-1}})=\chi(\mathbb{R}P^{n-1})=1.
    \end{equation*}
    If $m\not = 1$, choose any prime $p$ which divides $m$. Set $\bar{\tau}:= \bar{\sigma}^{m/p}$. For the action induced by $\bar{\tau}$ on the quotient space mentioned above, we must have
    \begin{equation*}
     0=\lambda(\bar{\tau})\equiv \lambda(\bar{\tau}^{p}) \quad \mod p.
    \end{equation*}
    Here, we have used the $\mod p$ Lefschetz Theorem (see \cite{dugundji}, Chap. 9, Thm. 3.4, p. 231). Thus, we arrive at the contradiction that
    \begin{equation*}
                  0\equiv 1 =\chi(\mathbb{R}P^{n-1})\quad \mod p .                                                                                                                                                                                                                         \end{equation*}
  Hence, the action induced by $\bar{\sigma}$ itself fixes some point. This fixed point corresponds precisely to a pencil $\Gamma$ of geodesics as asserted.
\section{Proof of Theorem \ref{thm-auto}.}
\begin{proof}
\indent

     Select an arbitrary element $\sigma \in Aut(\Omega_{k})$ having a prime order $p$. By \cite{krantz} (see Thm. 1.4.3, p. 22), we have
    \begin{equation}\label{eq-fix}
     \#Fix(\sigma) \leq 2 .
    \end{equation}
      Since $\Omega :=\Omega_{k}\backslash Fix(\sigma)$ is not compact, we shall work with \textit{generalised} Lefschetz numbers (as in \cite{dugundji}, Def. 3.1, p. 421). Note that $\Omega$ is a \textit{Lefschetz space} (by \cite{dugundji}, Chap. 5, Thm. 4.1, p. 423).\\

      Now, the restriction of $\sigma$ to $\Omega$ is an automorphism without fixed points. Therefore, the generalised Lefschetz number $\lambda(\sigma)$ of $\sigma$ must vanish. As $\Omega$ has the homotopy type of a wedge of $(k+ \#Fix(\sigma))$ circles, we obtain (using \cite{dugundji}, Chap. 9, Prop. 4.7. p. 233)
\begin{equation*}
 \chi(\Omega)=1-k-\#Fix(\sigma)
\equiv 0 \mod p.\end{equation*}
In view of Eq. (\ref{eq-fix}), this implies that $p$ must divide either $k$ or one of $k\pm 1$. Thus, the assertion of the Theorem follows.
\end{proof}
\begin{acknowledgement}
    The author would like to thank Chennai Mathematical Institute for support by a post-doctoral fellowship.
\end{acknowledgement}
\bibliographystyle{plain}
\bibliography{cartan's}
\addcontentsline{toc}{chapter}{Bibliography}
\end{document}